\documentclass[11pt]{amsart}

\usepackage{amssymb}

\usepackage{hyperref} %bib refs = green boxes, thm refs = red boxes

 \newtheorem{thm}{Theorem}[section]
  
 \newtheorem{cor}[thm]{Corollary}
 
 \newtheorem{prop}[thm]{Proposition}

 \theoremstyle{definition}
 \newtheorem{defn}[thm]{Definition}
 \newtheorem{rem}[thm]{Remark}
 
 \numberwithin{equation}{section}

\newcommand{\cref}[1]{Corollary~\ref{#1}}   %use: \cref{labelname}
\newcommand{\R}{\mathbb R}

\newcommand{\orb}{\mathcal{O}}
\newcommand{\mO}{\mathcal{O}}
\newcommand{\mX}{\mathcal{X}}
\newcommand{\mY}{\mathcal{Y}}

\newcommand{\SU}{\operatorname{SU}}
\newcommand{\SO}{\operatorname{SO}}

\newcommand{\Ind}{\ensuremath{\operatorname{Ind}}}
\newcommand{\Res}{\ensuremath{\operatorname{Res}}}

\newcommand{\Vol}{\operatorname{vol}}

\newcommand{\Spec}{\operatorname{Spec}}

\newcommand{\bs}{\backslash}

\begin{document}

\title{Equivariant isospectrality and Sunada's Method}
\author[C. J. Sutton]{Craig J. Sutton}

\address{
Dartmouth College\\
Department of Mathematics\\
Hanover, NH\\
USA}

\email{craig.j.sutton@dartmouth.edu}

\thanks{This work was partially supported by an NSF Postdoctoral Research Fellowship and NSF grant DMS 0605247}



\date{May 8, 2009}

\begin{abstract}
We construct pairs and continuous families of isospectral yet locally non-isometric 
orbifolds via an equivariant version of Sunada's method. 
We also observe that if a good orbifold $\mathcal{O}$ and a smooth manifold $M$ are isospectral, 
then they cannot admit non-trivial finite Riemannian covers $M_1 \to \mathcal{O}$ and $M_2 \to M$
where $M_1$ and $M_2$ are isospectral manifolds.
\end{abstract}

\maketitle

\section{Introduction}

An orbifold is a generalization of a manifold where
we allow the coordinate charts to be modeled on quotients of euclidean space.
To be more precise, an \emph{$n$-dimensional orbifold chart},
on a second countable Hausdorff space $\mO$, is a triple $(U, \Gamma_U, \pi_U)$, 
where $U$ is an open subset of $\mO$,
$\Gamma_U$ is a finite group of diffeomorphisms of $\R^n$,
and $\pi_U : \R^n \to U$ is a $\Gamma_U$-invariant mapping that  
induces a homeomorphism between the quotient space $ \Gamma_U \bs \R^n$ and $U$. 
An orbifold structure on $\mO$ is a collection $\{(U_\alpha, \Gamma_\alpha, \pi_\alpha)\}_{\alpha \in J}$
of $n$-dimensional orbifold charts where $\mO = \cup_{\alpha \in J} U_\alpha$ and the charts satisfy 
a compatibility condition.
A consequence of the compatibility condition is that it ensures that if $x$ is in the intersection of two charts 
$(U_1, \Gamma_1, \pi_1)$ and $(U_2, \Gamma_2, \pi_2)$, then for any 
$\tilde{x}_1 \in \pi_1^{-1}(x)$ and $\tilde{x}_2 \in \pi_2^{-1}(x)$, the isotropy groups 
$\Gamma_{1\tilde{x}_1}$ and $\Gamma_{2\tilde{x}_2}$ are isomorphic.
This common group is called the isotropy of $x$ and 
we will say that a point is \emph{singular} if it has non-trivial isotropy.

Orbifolds arise quite naturally in the context of group actions. 
Indeed, if $\Gamma$ is a group of diffeomorphisms acting properly discontinuously 
on a manifold $M$ with a fixed-point set of codimension at least two, 
then the quotient space $\Gamma \bs M$ is an orbifold \cite{Thurston}.
It is common to refer to an orbifold arising as a quotient of a manifold as a 
\emph{good} or \emph{global} orbifold, otherwise the orbifold is said to be \emph{bad}.

In the case of a good orbifold it is clear how one should define 
various analytic and geometric concepts. For instance, if $\orb =  \Gamma \bs M$
 is a good orbifold, then a function $f: \orb \to \R$ is said to be \emph{smooth} 
 if its lift to $M$ is a smooth function. 
 In this way we see that the space of smooth functions on $\orb$ 
 is naturally identified with the space of $\Gamma$-invariant smooth functions on $M$. 
 Since an arbitrary orbifold is locally a good orbifold, we will agree to say that $f: \orb \to \R$ 
 is smooth if its pullback on each local coordinate system is smooth. 
 Similarly, a Riemannian structure on a good orbifold $\orb = \Gamma \bs M$
 is defined by a $\Gamma$-invariant Riemannian metric on $M$ 
 (for a detailed discussion of the general case see \cite[Sec. 2]{Weilandt}).
 Continuing in this fashion we may extend the study of 
 geometric analysis to orbifolds. In particular, the classical inverse spectral
  problem carries over naturally to this setting.

Indeed, the Laplace operator $\Delta_\orb$ associated to a 
Riemannian orbifold $\orb$ will be defined locally through the 
coordinate charts, and as in the case of a Riemannian manifold one can see that 
the eigenvalues form a non-decreasing sequence of non-negative real numbers 
tending towards infinity \cite{Chiang}. This sequence is known as the \emph{spectrum} 
of the orbifold and as usual we will agree to say that two orbifolds sharing the same 
spectrum are \emph{isospectral}. Recently, it has become of interest to explore the 
relationship between the geometry of an orbifold and its spectrum.
For a nice introduction to orbifolds and the isospectral problem
we refer the reader to \cite{Weilandt}. 

This note is motivated by the question of whether one can ``hear'' 
the presence of a singularity.
That is, we wonder whether it is possible to construct a pair 
of isospectral orbifolds  where one orbifold has singular points  
while the other does not.
In support of the existence of such a pair, we note that there are examples of 
isospectral good orbifolds with a common cover for which the 
size of the maximal isotropy group differs \cite{RSW}. 
However, in \cite{GorRos} it was shown that whenever two isospectral \emph{good} 
orbifolds share a common Riemannian cover their respective singular sets are 
either both trivial or both non-trivial. Thus, setting aside bad orbifolds, 
one needs to examine isospectral orbifolds without common Riemannian coverings.
With this in mind we recast Sunada's method in the context of equivariant isospectrality.

Given a compact Lie group $G$, we will say that two 
Riemannian $G$-manifolds $M_1$ and $M_2$ are 
\emph{equivariantly isospectral} with respect to $G$
if there is a unitary isomorphism $U : L^{2}(M_{1}) \to L^{2}(M_{2})$ 
intertwining the Laplacians which is also an equivalence of the 
natural $G$-representation $(g\cdot f)(x) = f(g^{-1}\cdot x)$, where 
$g \in G$ and $f \in L^2(M_i)$ ($i = 1,2$). 
We then have the following equivariant version of Sunada's theorem.

\begin{thm}\label{Thm:MainSunada}
Let $(G, \Gamma_1, \Gamma_2)$ be a triple of finite groups where 
$\Gamma_1 , \Gamma_2 \leq G$ are subgroups 
such that for any conjugacy class $C \subset G$ we have 
$\#(\Gamma_1 \cap C ) = \# (\Gamma_2 \cap C)$.
Then, if $M_1$ and $M_2$ are $G$-equivariantly isospectral manifolds,
the orbifolds $\Gamma \bs M_1$ and $\Gamma_2 \bs M_2$ are isospectral.
(See Theorem~\ref{Thm:EquivSunada} for a more general statement.) 
\end{thm} 

Using the theorem above we obtain pairs and continuous families of isospectral 
good orbifolds without common Riemannian covers (see Section~\ref{Sec:Examples}).
These spaces arise as finite quotients of the 
equivariantly isospectral simply-connected spaces constructed in 
\cite{Gordon2} \cite{Schueth4} and \cite{Sut}, 
and appear to be candidates for demonstrating that 
you cannot hear the presence of a singularity.
However, the following theorem shows that singularities are ``audible'' 
within the class of isospectral orbifolds formed via Theorem~\ref{Thm:MainSunada}.

\begin{thm}\label{Thm:HearSing}
Let $\mO_1$ and $\mO_2$ be isospectral good orbifolds.
If $ \pi_1: M_1 \to \mO_1$ and $\pi_2: M_2 \to \mO_2$ are non-trivial 
finite Riemannian orbifold covers with $M_1$ and $M_2$ isospectral manifolds,
then $\mO_1$ has a singular point if and only if $\mO_2$ has a singular point.
\end{thm}

%%%%%%%%%%%%%%%%%%%%%%%%%%%%%%%%%%%%%%%%%%%%%%%
%%%%%%  Equivariant isospectrality and Sunada's method   %%%%%%%%%%%%%%
%%%%%%%%%%%%%%%%%%%%%%%%%%%%%%%%%%%%%%%%%%%%%%%

\section{An equivariant Sunada method}\label{sec:EquivSunada}

Let $G$ be a compact Lie group which acts via isometries on a Riemannian manifold $(M,g)$. 
Then $G$ has a natural representation on $L^{2}(M)$ where for each 
$f \in L^2(M)$ the function $g.f$ is given by 
$$(g.f)(x) \equiv f(g^{-1}x).$$ 
We denote this representation by $\tau^{G}$\label{RegularRep}. Since $G$ acts via isometries on $(M,g)$,
we see that $\tau^{G}$ commutes with the Laplacian $\Delta$. 
Hence, the decomposition of $L^2(M)$ into $\Delta$-eigenspaces, 
given by $L^{2}(M) = \oplus_{\lambda \in \Spec(\Delta)} L^{2}(M)_{\lambda}$, is invariant under $\tau^{G}$. 
Letting $\widehat{G}$ denote the set of equivalence classes of irreducible representations 
of $G$ we may also decompose $L^{2}(M)$ into $G$-invariant subspaces as follows 
$$L^{2}(M) = \oplus_{[\rho] \in \widehat{G}} L^{2}_{\rho}(M),$$ where $L^{2}_{\rho}(M)$ 
is the closed linear span of all irreducible subspaces of $L^{2}(M)$ on which $\tau^{G}$ 
is equivalent to $\rho$. Using the invariance of the $\Delta$-eigenspaces under $\tau^{G}$ 
we see that $\Delta$ preserves this decomposition and we let 
$\Spec(\Delta |_{L^{2}_{\rho}(M)}) = \{ \lambda^{\rho}_{1} \leq \lambda^{\rho}_{2} \leq \cdots \}$ 
denote the (possibly finite) spectrum of $\Delta : L^{2}_{\rho}(M) \to L^{2}_{\rho}(M)$.

\begin{defn}
Let $G$ be a compact Lie group acting via isometries on the Riemannian manifolds 
$(M_{1}, g_{1})$ and $(M_{2}, g_{2})$. Then $(M_{1}, g_{1})$ and $(M_{2}, g_{2})$ 
are said to be \emph{equivariantly isospectral} (with respect to $G$) if 
$\Spec(\Delta |_{L^{2}_{\rho}(M_{1})}) = \Spec(\Delta |_{L^{2}_{\rho}(M_{2})})$ 
for each $[\rho] \in \widehat{G}$. 
\end{defn}

\noindent
Equivalently, we have the following.

\begin{defn}
Two Riemannian $G$-manifolds $(M_{1}, g_{1})$ and $(M_{2}, g_{2})$ are said to 
be \emph{equivariantly isospectral} (with respect to $G$) if there exists a unitary map 
$U : L^{2}(M_{1}) \to L^{2}(M_{2})$ such that 
\begin{enumerate}
\item $U \circ \Delta_{1} = \Delta_{2} \circ U$; that is, $M_{1}$ and $M_{2}$ are isospectral.
\item $U \circ \tau_{1}^{G} = \tau_{2}^{G} \circ U$; that is, the natural representations are equivalent via $U$.
\end{enumerate}
\end{defn}

\begin{rem}
One can check that the isospectral spaces discussed in \cite{Gordon2, Schueth2, Schueth3, Sut}
and \cite[Theorem 4.1]{Schueth4} are equivariantly isospectral. 
\end{rem}

Before we state the equivariant Sunada method we need to introduce Pesce's 
notion of $K$-equivalence \cite{Pesce}. We begin by recalling that given two 
representations $\tau_{1} : G \to GL(V_{1})$ and $\tau_{2} : G \to GL(V_{2})$ 
the \emph{multiplicity} of $\tau_{1}$ in $\tau_{2}$ , denoted $[\tau_{2} : \tau_{1}]$, 
is defined to be $\dim ( \operatorname{Hom}_{G}(V_{1}, V_{2}) )$, where 
$\operatorname{Hom}_{G}(V_{1}, V_{2})$ is the set of bounded linear maps 
$T: V_{1} \to V_{2}$ such that $T \circ \pi_{1}(g) = \pi_{2}(g) \circ T$ for any $g \in G$. 
 Now,  let $G$ be a compact Lie group and let $\widehat{G}$ denote the set of 
 (equivalence classes of) irreducible representations of $G$. For any $K \leq G$ 
 closed we define $$\widehat{G}_{K} = \{ \rho \in \widehat{G} : [\Res^{G}_{K}(\rho) : 1_{K} ] \neq 0 \},$$
where $\Res^{G}_{K}(\rho)$ denotes the restriction of $\rho$ to $K$ and $1_{K}$ 
denotes the trivial representation of $K$. Therefore, $\widehat{G}_{K}$ is the set 
of irreducible representations of $G$ which have non-trivial $K$-fixed vectors. 

\begin{defn}
Let $(\tau_{1}, V_{1})$ and $(\tau_{2}, V_{2})$ be two representations of $G$ 
such that $[\tau_{i} : \rho] < \infty$ for any $\rho \in \widehat{G}$ ($i=1,2$) and 
let $K \leq G$ be a closed subgroup. We will say that $\tau_{1}$ and $\tau_{2}$ 
are \emph{$K$-equivalent} representations, denoted $\tau_{1} \sim_{K} \tau_{2}$, 
if $[\tau_{1} : \rho] = [\tau_{2}: \rho]$ for each $\rho \in \widehat{G}_{K}$; that is, 
the restrictions of $\tau_{1}$ and $\tau_{2}$ to the smallest $G$-invariant subspaces 
of $V_{1}$ and $V_{2}$ (respectively) which contain all of the $K$-fixed vectors are equivalent.
In the case where $K$ is trivial we obtain the usual notion of equivalence and write
$\tau_1 \sim \tau_2$.
\end{defn}

\begin{defn}
Let $G$ be a compact Lie group.
\begin{enumerate}
\item Given a closed subgroup $H \leq G$ we define the \emph{quasi-regular} 
representation of $G$ on $L^{2}(G/H)$, denoted $\pi_{H}^{G}$,  via 
$$(\pi_{H}^{G}(g)f)(x) \equiv f(g^{-1} \cdot x),$$ 
for any $g \in G$ and $f \in L^{2}(G/H)$. 
That is, $\pi_{H}^{G}$ is the representation of $G$ induced by the trivial representation of $H$, 
which is often denoted by $\Ind^{G}_{H}(1_{H})$. 
We refer the reader to \cite[Sec. 6.1]{Folland} for a discussion of induction.
\item Given subgroups $K, H_{1}, H_{2} \leq G$ we will say that $H_{1}$ and $H_{2}$ 
are \emph{$K$-equivalent} subgroups if and only if $\pi_{H_{1}}^{G} \sim_{K} \pi_{H_{2}}^{G}$. 
In the case where $K$ is trivial we say $H_{1}$ and $H_{2}$ are \emph{representation equivalent} 
subgroups.

\item If $G$ is a finite group, then subgroups $H_1, H_2 \leq G$ are 
said to be \emph{almost conjugate} or \emph{Gassmann-Sunada equivalent} 
if and only if $\#([g]_G \cap H_{1}) = \#([g]_G \cap H_{2})$ 
for any $g \in G$, where $[g]_G$ denotes the conjugacy class of $g$ in 
$G$ (see \cite{Sunada}). 
In this case we call $(G, H_1, H_2)$ a \emph{Gassmann-Sunada triple}.
\end{enumerate}
\end{defn}
 
 \begin{rem}\label{Rem:GassmannSunada}
 If $G$ is finite, then $(G, H_1, H_2)$ is a Gassmann-Sunada triple if 
 and only if $\pi_{H_{1}}^{G} \sim \pi_{H_{2}}^{G}$.
 \end{rem}
 
We now state and prove our equivariant Sunada technique.

\begin{thm}\label{Thm:EquivSunada}
Let $M_{1}$ and $M_{2}$ be two (possibly isometric) isospectral 
Riemannian manifolds and let $G$ be a compact Lie group such that

\begin{enumerate}
\item $G$ acts by isometries on $M_{1}$ and $M_{2}$.
\item $M_{1}$ and $M_{2}$ are equivariantly isospectral with respect to $G$.
\item The actions of $G$ on $M_{1}$ and $M_{2}$ have the same generic stabilizer $K \leq G$.
\end{enumerate}
Now suppose that $H_{1}, H_{2} \leq G$ are closed, $K$-equivalent 
subgroups which act on $M_{1}$ and $M_{2}$ respectively such that 
the Riemannian submersions  $$\pi_{i} : M_{i} \to H_i \bs M_{i} \; (i=1,2),$$ 
have minimal fibers, where $H_i \bs M_i$ has the induced metric.
Then $H_1 \bs M_{1}$ and $H_2 \bs M_{2}$ 
are isospectral on functions. 
\end{thm}

\begin{rem}
In the case where $M_{1}$ and $M_{2}$ are isometric and $G$ is a finite group 
which acts freely we obtain Sunada's method. 
\end{rem}

\begin{rem}
The proof is essentially the same as in \cite{Sut}, where we demonstrated 
that one can construct simply-connected, \emph{locally} non-isometric isospectral 
spaces through a generalization of Sunada's method. 
\end{rem}

\begin{proof}
For each $i = 1,2$ and $\lambda \in \Spec (M_{i})$ we let $\tau_{i, \lambda}^{G}$ and 
$\tau_{i, \lambda}^{H_{i}}$ denote the natural representations of $G$ and $H_{i}$ on 
$L^{2}(M_{i})_{\lambda}$ respectively (see p.~\pageref{RegularRep}). Since the 
Riemannian submersions $\pi_{1} : M_{1} \to N_{1}$ and $\pi_{2} : M_{2} \to N_{2}$ 
have minimal fibers it follows that  for each $\lambda \in \Spec(M_{1}) = \Spec(M_{2})$ 
we have  $$\dim L^{2}(H_i \bs M_{i})_{\lambda} = [ \tau_{i, \lambda}^{H_{i}} : 1_{H_{i}}]  \;\; (i = 1,2).$$
Hence, $H_1 \bs M_{1}$ and $H_2 \bs M_{2}$ are isospectral if and only if 
$[ \tau_{1, \lambda}^{H_{1}} : 1_{H_{1}}] = [ \tau_{2, \lambda}^{H_{2}} : 1_{H_{2}}]$ 
for each $\lambda \in \Spec(M_{1}) = \Spec(M_{2})$.  Now, since 
$\tau_{i, \lambda}^{H_{i}} = \Res_{H_{i}}^{G} ( \tau_{i, \lambda}^{G})$, where $\Res^{G}_{H}(\rho)$ 
denotes the restriction of the representation $\rho$ of $G$ to the closed subgroup $H \leq G$, 
we obtain using Frobenius' reciprocity theorem:

\begin{eqnarray*}
[ \tau_{i, \lambda}^{H_{i}} : 1_{H_{i}}]	&=& [\Res^{G}_{H_{i}}(\tau^{G}_{i, \lambda}) : 1_{H_{i}}] \\
					&=& [\Res^{G}_{H_{i}}(\sum_{\rho \in \widehat{G}}[\tau^{G}_{i, \lambda}:\rho]\rho) : 1_{H_{i}}] \\
					&=& \sum_{\rho \in \widehat{G}}[\tau^{G}_{i, \lambda}:\rho][\Res^{G}_{H_{i}}(\rho) : 1_{H_{i}}] \\
					&=& \sum_{\rho \in \widehat{G}}[\tau^{G}_{i, \lambda}:\rho][\Ind^{G}_{H_{i}}(1_{H_{i}}) : \rho] \\
					&=& \sum_{\rho \in \widehat{G}}[\tau^{G}_{i, \lambda}:\rho][\pi^{G}_{H_{i}} : \rho]. \\
\end{eqnarray*}

We now recall the following theorem of Donnelly.

\begin{thm}[\cite{Donn}, p. 25]
Let $G$ be a compact Lie group and $X$ a compact, smooth $G$-space 
with principal orbit type $G/K$; that is, $K$ is the 
generic stabilizer of the $G$-action on $X$. Then the decomposition of 
$L^{2}(X)$ into $G$-irreducibles contains precisely those finite dimensional 
representations appearing in the decomposition of $\pi^{G}_{K} = \Ind^{G}_{K}(1_{K})$ 
the quasi-regular representation of $G$ with respect to $K$. Also, if the orbit 
space $G \bs X$ has dimension greater than 1, then each irreducible 
appears an infinite number of times.
\end{thm}

Using Frobenius' theorem once again we have $[\pi^{G}_{K}: \rho] = [ \Res^{G}_{K}(\rho) : 1_{K}]$ 
for each $\rho \in \widehat{G}$; hence, from Donnelly's result we conclude that for each $i = 1,2$ we have 
$$[\tau^{H_{i}}_{i, \lambda} : 1_{H_{i}}] = \sum_{\rho \in \widehat{G}_{K}}[\tau^{G}_{i, \lambda}:\rho][\pi^{G}_{H_{i}} : \rho],$$
where $\widehat{G}_{K} = \{ \rho \in \widehat{G} : [\Res_{K}^{G}( \rho ) : 1_{K} ] \neq 0 \}$. 
Since $(M_{1}, g_{1})$ and $(M_{2}, g_{2})$ are equivariantly isospectral and 
$\Ind_{H_{1}}^{G}(1_{H_1})$ and $\Ind_{H_{2}}^{G}(1_{H_2})$ are $K$-equivalent representations, it follows that 
$[ \tau_{1, \lambda}^{H_{1}} : 1_{H_{1}}] = [ \tau_{2, \lambda}^{H_{2}} : 1_{H_{2}}]$ 
for all $\lambda \in \Spec(\Delta_{M_{1}}) = \Spec(\Delta_{M_{2}})$. 
Hence, $H_{1} \bs M_{1}$ and $H_{2} \bs M_{2}$ are isospectral.

\end{proof}

\begin{proof}[Proof of Theorem~\ref{Thm:MainSunada}]
By Remark~\ref{Rem:GassmannSunada} this is just a special case of Theorem~\ref{Thm:EquivSunada}.
\end{proof}

%%%%%%%%%%%%%%%%%%%%%%%%%%%%%%%%%%%%
%%%%%%%%   Sunada isospectral orbifolds %%%%%%%%%%%%
%%%%%%%%%%%%%%%%%%%%%%%%%%%%%%%%%%%%

\section{Sunada isospectral orbifolds with different local geometries}\label{Sec:Examples} 

In this section we use the equivariant Sunada method to construct examples
of isospectral good orbifolds without a common Riemannian cover. First, we recall 
that in \cite{Sut} we showed that there is a connected compact 
semisimple Lie group $H$ with the following properties (\cite[Cor. 3.4]{Sut}):

\begin{enumerate}
\item $H$ admits faithful representations $\rho_1 : H \to \SU(n)$ and 
$\rho_2: H \to \SU(n)$ for all $n$ greater than some $N_H$;
\item $H_1 = \rho_1(H)$ and $H_2 = \rho_2(H)$ are 
representation equivalent subgroups of $\SU(n)$ 
that are not conjugate via any automorphism;
\item the manifolds $\SU(n)/H_1$ and $\SU(n)/H_2$ are simply-connected. 
\end{enumerate} 

\noindent
It follows that the Riemannian manifolds $(\SU(n)/H_1, g_1)$ and $(\SU(n)/H_2, g_2)$, 
where $g_1$ and $g_2$ are the metrics induced by the 
bi-invariant metric on $\SU(n)$, are isospectral yet locally non-isometric \cite[Thm. 3.6]{Sut}.
In fact, one can readily see that these spaces are equivariantly isospectral with respect to $\SU(n)$. 
Hence, for any finite representation equivalent subgroups 
$\Gamma_1, \Gamma_2 \leq \SU(n)$ 
one obtains a pair of locally non-isometric isospectral orbifolds 
$\Gamma_1 \backslash \SU(n) / H_1$ and $\Gamma_2 \backslash \SU(n) / H_2$.
In particular, we may take $\Gamma_1 = \Gamma_2 \leq \SU(n)$.
Alternatively, we may proceed as follows.

Let $H$ be a Lie group as in the previous paragraph and let 
$(G, \Gamma_1, \Gamma_2)$ be a Gassmann-Sunada triple. 
Then for any $n \geq N_H$ that is sufficiently large we may find 
an injective homomorphism of $G$ into $\SU(n)$. 
For example, if $(V, \rho)$ is a faithful unitary 
representation of the finite group $G$, then 
$$g \mapsto \left(\begin{array}{cc}\rho(g) & 0 \\0 & \det (\rho(g)^{-1})\end{array}\right)$$
is a (non-trivial) homomorphism of $G$ into $\SU( \dim(V) + 1)$ 
and letting $n \geq \max(N_H, \dim(V) +1)$ we obtain the 
appropriate representation. It then follows that $(\SU(n)/H_1, g_1)$ 
and $(\SU(n)/H_2, g_2)$ are $G$-equivariantly isospectral and 
we conclude from Theorem~\ref{Thm:MainSunada} that 
$(\Gamma_1 \backslash \SU(n)/H_1, g_1)$ and 
$(\Gamma_2 \backslash \SU(n)/H_2, g_2)$ are isospectral quotient spaces.
If $n$ is taken to be sufficiently large we can arrange for 
both quotient spaces to be manifolds. 

We summarize the preceding discussion as follows.

\begin{prop}\label{Prop:GSExas}
Let $(G, \Gamma_1, \Gamma_2)$ be a Gassmann-Sunada triple and let $H$ be as above. 
Then there is an $N = N(G, H)$ such that for $n \geq N$, the normal homogeneous spaces
$SU(n)/H_1$ and $\SU(n)/H_2$ 
admit effective $G$-actions with respect to which the quotients 
$\Gamma_1 \bs \SU(n) / H_1$ and $\Gamma_2 \bs \SU(n)/H_2$ 
are isospectral yet locally non-isometric.
\end{prop}

\begin{rem}
In the above one can also consider any left-invariant metric $g$ on $\SU(n)$ 
which is also right-invariant with respect to $H_1$ and $H_2$. However, we do not know 
(except in the case where $g$ is bi-invariant) whether the resulting quotient spaces will be non-isometric.
It is also still unknown whether the underlying topological spaces $\SU(n)/H_1$ and $\SU(n)/H_2$ 
are homeomorphic. 
 If they prove to be non-homeomorphic, then the isospectral manifolds $(\SU(n)/H_1, g_1)$ 
 and $(\SU(n)/H_2, g_2)$, where $g_1$ and $g_2$ are normal homogeneous,
would demonstrate that the topological universal cover is not a spectral invariant. 
 All other known isospectral manifolds have homeomorphic universal covers.
\end{rem}

\begin{rem}
In \cite{Schueth4} pairs of isospectral metrics 
were constructed on the $5$-dimensional sphere. These metrics are 
equivariantly isospectral with respect to a non-free isometric action of the $2$-torus $T^2$.
It follows from Theorem~\ref{Thm:EquivSunada} that for any finite group $\Gamma \leq T^2$ one obtains 
isospectral metrics on the spherical orbifold
$\Gamma \bs S^5$ and for certain $\Gamma$ this orbifold has a non-trivial singular set.
\end{rem}

As we noted in the introduction, it is an interesting question to determine whether one can discern 
the presence of a singularity from knowledge of the spectrum.
The preceding collection of isospectral quotient spaces 
appears to hold some promise of containing an example of a manifold 
that is isospectral to a space with singularities. 
However, Theorem~\ref{Thm:HearSing}, the proof of which is given below, ensures that no such example exists 
among the pairs $\Gamma_1 \bs \SU(n)/H_1$ and $\Gamma_2 \bs \SU(n)/H_2$
constructed above.

Before we give the proof of Theorem~\ref{Thm:HearSing} 
we recall the definition of an orbifold covering map.

\begin{defn}
Let $\mX$ and $\mY$ be two $n$-dimensional Riemannian orbifolds. 
A mapping $p : \mX \to \mY$ is said to be a \emph{Riemanian orbifold 
covering} if for each $y \in \mY$ there is a coordinate chart 
$(U, \Gamma_U, \pi_U)$ 
containing $y$ such that: 
\begin{enumerate}
\item $p^{-1}(U)$ is a disjoint union of coordinate charts 
$\{(V_\alpha , \Gamma_\alpha, \pi_\alpha)\}_{\alpha \in J}$;
\item for each $\alpha \in J$ there is a monomorphism $i_\alpha : \Gamma_\alpha \to \Gamma_U$;
\item for each $\alpha \in J$ there is an isometry 
$\tilde{p}_{\alpha}: (\R^n, g_\alpha) \to (\R^n, g_U)$,
where $g_\alpha$ and $g_U$ are the Riemannian structures corresponding  
to the coordinate charts $V_\alpha$ and $U$ respectively, such that 
$\tilde{p}_\alpha(\gamma \cdot x) = i_\alpha(\gamma)\cdot \tilde{p}_\alpha(x)$
for each $\gamma \in \Gamma_\alpha$ and $x \in \R^n$,
and $\pi_U \circ \tilde{p}_\alpha = p_\alpha \circ \pi_{\alpha}$, where 
$p_\alpha$ is the restriction of $p$ to $V_\alpha$.
\end{enumerate}
\end{defn}

\begin{rem}
Good orbifolds are characterized by the property that they 
admit  (Riemannian) orbifold covers that are manifolds. 
For more details on orbifold coverings the reader is 
encouraged to consult \cite[Sec. 3]{Choi}.
\end{rem}

\begin{proof}[Proof of Theorem~\ref{Thm:HearSing}]
Given an orbifold $\mY$ we let $\Sigma_\mY$ denote 
the collection of singular points of $\mY$.
Then $\mY - \Sigma_\mY$ is an open and dense set.
Now suppose $p: \mX \to \mY$ is a Riemannian orbifold cover 
of complete spaces where each fiber is  countable.
(For instance, $p$ could be the quotient map associated to 
a properly discontinuous action of a discrete group on a manifold.)
Then for any $x_0 \in \pi^{-1}(\mY - \Sigma_\mY)$
we may define the \emph{Dirichlet fundamental domain of 
$p$ with center $x_0$} to be the set 
$C_{x_0} = \{ x \in \mX: d(x_0, x) < d(x_0, x') \mbox{ for any } x' \neq x \in F_x \}$, 
where $F_x$ denotes the fiber of $p$ through $x$. Fixing $x_0$ as above, 
one can check that the fundamental domain has the following properties (cf. \cite[Proposition 1.9.29]{Eber}):
\begin{enumerate} 
\item $\{C_x\}_{x \in F_{x_0} }$ is a collection of pairwise disjoint open sets;
\item $\mX = \cup_{x \in F_{x_0}} \overline{C}_x$;
\item The measure of the boundary of $C_x$ is zero for each $x \in F_{x_0}$;
\item $\Vol(C_x) = \Vol(\mO)$ for all $x \in F_{x_0}$.
\end{enumerate}

Now, let $\pi_1: M_1 \to \mathcal{O}_1$ and $\pi_2: M_2 \to \mathcal{O}_2$ 
be as in the hypotheses,
and for each $i = 1, 2$, let $C_i \subset M_i$ be the Dirichlet fundamental domain of 
$\pi_i$ centered at some $x_i \in M_i - \pi_i^{-1}(\Sigma_{\mO_i})$. Then we see that 
$\Vol(M_i) = d_i \cdot \Vol(C_i) = d_i \cdot \Vol(\mO_i)$, 
where $d_i$ is the order of $F_{x_i}$, 
and, using the fact that volume is a spectral invariant, we conclude from our 
hypotheses that $d_1 = d_2$.

We now recall that for any closed good orbifold $\mathcal{O}$ with eigenvalue 
spectrum  $\lambda_1 = 0 < \lambda_2 \leq \lambda_3 \leq \cdots \nearrow + \infty$,
we have the following asymptotic expansion of the heat trace due to Donnelly \cite{Donn2} (cf. \cite[Theorem 4.8]{DGGW}):
$$\sum_{i = 1}^{+\infty} e^{-\lambda_i t} \stackrel{t \searrow 0}{\sim} 
(4\pi t)^{-\frac{\dim (\mathcal{O})}{2}} \sum_{k = 0}^{+ \infty} a_k t^k + \sum_{S} B_S(t),$$
where $S$ varies over the strata of the singular set of $\mathcal{O}$ and where 
$$B_{S}(t) = (4\pi t)^{-\frac{\dim(S)}{2}} \sum_{k=0}^{+\infty} b_{k,S} t^k$$ with $b_{0,S} \neq 0$.
The coefficients $a_k$ in the first part of the 
asymptotic expansion above are given in terms of 
integrals (with respect to the Riemannian density) of expressions in the curvature of $\mO$ and its covariant derivatives.
Now, if $\pi: M \to \mO$ is a finite Riemannian cover of degree $d$ and
$\{\tilde{a}_k\}_{k \geq 0}$ denotes the corresponding heat invariants of $M$, then 
$\tilde{a}_k = d \cdot a_k$.
Therefore, if we let $a_{i,k}$ denote the corresponding terms of Donnelly's asymptotic expansion
for the orbifold $\mO_i$, 
it follows from the isospectrality of $M_1$ and $M_2$ (and the fact that 
$d_1 = d_2$) 
that the heat invariants $a_{1,k}$ and  $a_{2,k}$ are equal for each non-negative integer $k$. 
It then follows,  since $b_{0,S}$ is positive for each singular strata $S$ \cite[p. 218-220]{DGGW}, 
that $\mO_1$ has singular points if and only if $\mO_2$ has singular points.
\end{proof}

\begin{rem}
In \cite{Donn2} there is a typographical error in the statement of Donnelly's asymptotic expansion.
In the proof above, we have used the correct statement as found in \cite{DGGW} 
where the asymptotic expansion is generalized to include \emph{all} orbifolds.
\end{rem}

\begin{cor}\label{Cor:HearSing}
Let $M_1$ and $M_2$ be two isospectral closed Riemannian manifolds and let $\Gamma_1$ and $\Gamma_2$ be two 
discrete groups acting properly discontinuously and isometrically on $M_1$ and $M_2$ respectively.
If the quotient spaces $\Gamma_1 \bs M_1$ and $\Gamma_2 \bs M_2$ are isospectral, then 
$\Gamma_1$ acts freely on $M_1$ if and only if $\Gamma_2$ acts freely on $M_2$.
\end{cor}

We note that Donnelly's asymptotic expansion of the heat trace also demonstrates  
that if a good orbifold $\mathcal{O}$ with singular points and a manifold $M$ have a common 
Riemannian covering, then they cannot be isospectral  \cite{GorRos}. 
It is also shown in \cite{DGGW}, through an asymptotic expansion
of the heat trace
valid for \emph{all} orbifolds, that if an even (respectively, odd) 
dimensional orbifold has a singular strata of odd (respectively, even) dimension, 
then it cannot be isospectral to a smooth manifold. 

We conclude this note with the following observation concerning 
isospectral metrics on spherical orbifolds.

\begin{prop}\label{Prop:Spherical}
For each $n \geq 8$ there are spherical orbifolds of dimension $n$ 
that admit multiparameter families of isospectral yet locally non-isometric metrics.
\end{prop}

\begin{proof}
We recall the following method due to Gordon.

\begin{thm}[ \cite{Gordon2} Thm. 1.2]\label{Thm:TorusMethod}
Let $T$ be a torus. Suppose $T$ acts by isometries on two compact Riemannian 
manifolds $M_{1}$ and $M_{2}$ and that the action of $T$ on the principal orbits is free. 
Let $M_{i}'$ be the union of all the principal orbits in $M_{i}$, so $M_{i}'$ is an open 
submanifold of $M_{i}$ and a principal $T$-bundle, $i = 1,2$. For each subtorus $K \leq T$ 
of codimension at most one, suppose that there exists a diffeomorphism $\tau_{K} : M_{1} \to M_{2}$ 
which intertwines the actions of $T$ and which induces an isometry $\bar{\tau}_{K}$ 
between the induced metrics on the quotient manifolds $K \backslash M_{1}'$ and 
$K \backslash M_{2}'$. Assume further that the isometry $\bar{\tau}_{K}$ satisfies 
$\bar{\tau}_{K*}(\overline{H}_{K}^{1}) = \overline{H}_{K}^{2}$, where $\overline{H}_{K}^{i}$ 
is the projected mean curvature vector field for the submersion $M_{i}' \to K \backslash M_{i}'$. 
Then in the case that $M_{1}$ and $M_{2}$ are closed, they are isospectral. In the case where 
$M_{1}$ and $M_{2}$ have boundary, then they are Dirichlet isospectral, and under the 
additional assumption that $\partial(M_{i} )\cap M_{i}'$ is dense in $\partial(M_{i}')$ 
($i=1,2$), the manifolds are also Neumann isospectral. 
\end{thm}

In the proof of the above Gordon constructs an explicit intertwining operator 
$Q: L^2(M_1) \to L^2(M_2)$ of the Laplacians $\Delta_1$ and $\Delta_2$ which 
one can see is $T$-equivariant. Hence, the isospectral manifolds constructed 
via Thorem~\ref{Thm:TorusMethod} are $T$-equivariantly isospectral.
 Gordon then considers for $n \geq 8$ the standard $n$-sphere 
$S^n \subset \R^{n+1} = \R^2 \oplus \R^2 \oplus \R^{n-3}$
with the natural action of $T^2 = \SO(2) \oplus \SO(2) \oplus I_{n-3}$, 
and uses Theorem~\ref{Thm:TorusMethod} to construct multiparameter families 
of locally non-isomsetric $T^2$-equivariantly isospectral metrics on $S^n$ \cite[Cor. 3.10]{Gordon2}. 
Since each element of $T^2$ fixes at least a sphere of dimension $n-4$, 
we see that the $T^2$-action is not free and it follows from our 
equivariant Sunada method (Theorem~\ref{Thm:EquivSunada})
that for any finite subgroup $\Gamma \leq T^2$ the good orbifold $\Gamma \bs S^n$ admits a 
non-trivial multiparameter family of locally non-isometric metrics.
\end{proof}

\begin{rem}
We note that instead of using the equivariant Sunada technique, 
the above can also be seen by applying a simple perturbation argument to 
Gordon's sphere examples.
\end{rem}

\subsection*{Acknowledgments}
We thank the referee for suggesting that we consider Gassmann-Sunada triples in Section~\ref{Sec:Examples}.

%%%%%%%%%%%%%%%%%%%%%%%%%%%%%
%%%%%%%%%% Bibliography  %%%%%%%%%%%
%%%%%%%%%%%%%%%%%%%%%%%%%%%%%

\bibliographystyle{amsalpha}

\end{document}